\documentclass[amssymb,11pt]{amsart}
\usepackage{amscd,amssymb,euscript,amsthm}
\theoremstyle{plain}
\newtheorem{thm}{Theorem}[section] 

\newtheorem{prop}[thm]{Proposition}
\newtheorem{lem}[thm]{Lemma}
\theoremstyle{definition}
\newtheorem{defn}[thm]{Definition}
\theoremstyle{remark}
\newtheorem{rem}[thm]{Remark}

\numberwithin{equation}{section}
\newcommand{\id}{\operatorname{id}}

\newcommand{\ran}{\operatorname{ran}}

\def\<{\left<}
\def\>{\right>}
\def\cstar{$C^*$-algebra}
\begin{document}
\title[Asymptotic Stability]{Asymptotic Stability I:\\
Completely positive maps}
\author{William Arveson}
\thanks{supported by 
NSF grant DMS-0100487} 
\address{Department of Mathematics,
University of California, Berkeley, CA 94720}
\email{arveson@math.berkeley.edu}
%
\begin{abstract} 
We show 
that for every ``locally finite" 
unit-preserving completely positive map $P$ acting 
on a \cstar, 
there is a corresponding $*$-automorphism $\alpha$ of 
another unital \cstar\ such that the  two 
sequences $P, P^2,P^3,\dots$ and 
$\alpha, \alpha^2,\alpha^3,\dots$ have the same 
{\em asymptotic} behavior.  The automorphism   
$\alpha$ is uniquely determined by $P$ up to conjugacy.  
Similar results hold for normal completely positive 
maps on von Neumann algebras, as well as for 
one-parameter semigroups.  

These results are operator algebraic counterparts 
of the classical theory of Perron and Frobenius on 
the structure of 
square matrices with nonnegative entries.  
\end{abstract}
\maketitle

\section{Introduction}\label{S:intro}

The purpose of this paper is to show 
that many completely positive maps on \cstar s, and  
normal completely positive maps on von Neumann algebras,  
have significant 
asymptotic stability properties.  
This material arose in connection with 
our work on an asymptotic spectral invariant 
for single automorphisms of 
\cstar s, and 
for one-parameter semigroups of endomorphisms of von Neumann 
algebras.  Those applications will be taken up 
in a subsequent paper.  However, since they provide 
the motivation and conceptual foundation 
for the discussion below, we offer the 
following remarks concerning the problem that inspired 
this work, and describe the connections 
between the results of this paper and noncommutative dynamics.  

Let $n\geq 2$ be a positive integer.  The noncommutative 
Bernoulli shift of rank $n$ is the automorphism of the 
hyperfinite $II_1$ factor $R$ that is associated 
with the bilateral shift acting on the UHF algebra
$$
\mathcal A_n=\bigotimes_{k=-\infty}^{+\infty} A_k,  
$$
each factor $A_k$ being the \cstar\ of $n\times n$ 
matrices.  The GNS construction applied to the 
tracial state gives rise 
to a representation of $\mathcal A_n$ whose weak closure 
is $R$; and since the shift on $\mathcal A_n$ 
preserves the trace  
it can be extended naturally to a 
$*$-automorphism of $R$,  which we 
denote by $\sigma_n$.  The problem of whether $\sigma_m$ 
is conjugate to $\sigma_n$ for $m\neq n$ was 
solved by Connes and St\o rmer \cite{conSt} by introducing  
a  noncommutative generalization of 
the Kolmogorov-Sinai entropy of ergodic 
theory.  The entropy of the noncommutative shift of 
rank $n$ was computed in \cite{conSt}, and was found 
to have the value $\log n$, thereby settling 
the issue of conjugacy of the various shifts 
$\sigma_n$.  

We now consider another construction 
of automorphisms of von Neumann algebras that 
resembles the construction of 
finite state Markov processes 
from their transition probability matrices.  
This construction 
begins with a pair $(A,P)$, consisting of a finite 
dimensional \cstar\ $A$ and a completely positive 
linear map $P:A\to A$ satisfying $P(\mathbf 1)=\mathbf 1$.  
For technical 
reasons we require that 
$P(e)\neq 0$ for every minimal central projection 
$e\in A$; there is no essential loss if one thinks of $A$ as 
a full matrix algebra $M_n(\mathbb C)$ - and in that case 
the technical hypothesis is automatically satisfied.  There is a 
``noncommutative Markov process" that can be constructed 
from  $(A,P)$ as follows.  Briefly, a 
dilation theorem of Bhat \cite{bhatIndex, bhatMin}, as formulated 
in Chapter 8 of \cite{arvMono}, gives 
rise to a pair $(M_0, \sigma_0)$ consisting of a von Neumann 
algebra $M_0$ and a normal $*$-endomorphism $\sigma_0:M_0\to M_0$ 
that is appropriately related to the pair $(A,P)$.  The 
technical hypothesis implies that $\sigma_0$ is isometric, 
and one may then show that the endomorphism $\sigma_0$ can 
be extended appropriately to a $*$-automorphism 
of a larger von Neumann algebra $M\supseteq M_0$.  Let us 
denote the latter automorphism by $\sigma^P$.  

By analogy with the theory of noncommutative Bernoulli
shifts, we were led to conjecture that $\sigma^P$ and 
$\sigma^Q$ were generically not conjugate.  
The proof of that called for a new invariant, 
since the Connes-St\o rmer entropy 
is inappropriate for two reasons.  First, $\sigma^P$ 
typically acts on a type $I$ von Neumann algebra $M$, 
and second,  the construction of $M$ 
involves free products of copies of $A$ 
and not tensor products \cite{arvNCgen}, 
(see \cite{storSurv} for the significance of that fact).  
There is an asymptotic invariant for automorphisms $\alpha$ 
of von Neumann algebras (as well as for \cstar s) that 
we call the {\em asymptotic spectrum} $Sp_\infty(\alpha)$.  
We do not define $Sp_\infty(\alpha)$ here, but we do 
point out that it is a subset of the unit circle, perhaps 
finite, and is typically not closed.  The fact is 
that the asymptotic spectrum is relatively easy 
to compute and serves to distinguish 
between the various $\sigma^P$.  The key result 
on the computation of 
the asymptotic spectrum  is 
the following.

\begin{thm}\label{thm00}Let $(A,P)$ be a pair consisting of a 
von Neumann algebra $A$ and a normal completely positive 
map $P: A\to A$ satisfying $P(\mathbf 1)=\mathbf 1$, 
and let $(M,\sigma^P)$ be the associated $W^*$-dynamical 
system.  Then under appropriate 
hypotheses that include all pairs 
$(A,P)$ with finite-dimensional $A$, one has 
$$
SP_\infty(\sigma^P)=\sigma_p(P)\cap \mathbb T, 
$$
where $\sigma_p(P)$ denotes the point spectrum of $P$.  
\end{thm}

Thus, $\sigma^P$ and $\sigma^Q$ are not conjugate 
whenever $P$ and $Q$ have a different set of eigenvalues 
on the unit circle, and hence there is a continuum 
of non-conjugate automorphisms $\sigma^P$.  
The proof of Theorem \ref{thm00} has two components: a)
the development of properties of $Sp_\infty(\alpha)$, 
and b) an analysis of the asymptotic stability properties 
of completely positive maps.  In this paper we concentrate 
on b), the discussion a)  
will be taken up elsewhere.  

We now describe the contents of this paper in 
somewhat more detail.  We are concerned with the 
{\em asymptotic} behavior of the powers of a 
completely positive map $P: A\to A$, where $A$ is 
either a \cstar\ or a von Neumann algebra (in which 
case $P$ is assumed to be normal).  In order to illustrate 
the simplest case of that phenomenon, consider a normal 
completely positive map $P$ of 
$\mathcal B(H)$ into itself satisfying $P(\mathbf 1)=\mathbf 1$.  
In this setting there is a natural generalization of 
the ergodic-theoretic notion of mixing:  $P$ is said 
to be {\em mixing} if there is a normal state 
$\omega$ on $\mathcal B(H)$ with the property 
\begin{equation}\label{Eq01}
\lim_{n\to\infty}\omega(AP^n(B))=\omega(A)\omega(B),\qquad 
A,B\in \mathcal B(H).  
\end{equation}
If such an $\omega$ exists, then it must be invariant in 
the sense that $\omega\circ P=\omega$; indeed,  
$\omega$ is the {\em unique} normal $P$-invariant state.  
In case $\omega$ is faithful, it is not hard to see  
that (\ref{Eq01}) is equivalent to the following somewhat 
stronger absorption property: For every normal state 
$\rho$ on $\mathcal B(H)$ one has 
\begin{equation}\label{Eq02}
\lim_{n\to \infty}\|\rho\circ P^n-\omega\|=0.    
\end{equation}
The formula (\ref{Eq02}) represents the simplest form of 
the asymptotic stability that such maps $P$ can 
have.  Our objective is to develop a generalization 
of this kind of stability that is flexible enough to 
apply to a broad class of 
completely positive maps on von Neumann algebras 
(Theorem \ref{thm4}) and 
\cstar s (Theorem \ref{thm3}).

In order to keep the discussion as simple and 
focused as possible, we shall fix attention 
on pairs $(A,P)$ consisting 
of a \cstar\ $A$ (perhaps without unit) and a completely 
positive map $P: A\to A$ satisfying $\|P\|= 1$, making 
occasional comments about how formulations must be 
modified for normal maps on von Neumann algebras.  In 
this case, the idea of asymptotic stability is formulated 
as follows.  

The most rigid completely positive maps 
are called {\em quasiautomorphisms} below.  Roughly 
speaking, a quasiautomorphism is a completely 
positive contraction $Q: A\to A$ whose behavior away 
from its null space 
$\ker Q=\{z\in A: Q(z)=0\}$ is identical 
to that of a $*$-automorphism $\alpha$ of a secondary 
\cstar\ $B$ in the following sense: The powers of 
$\alpha$ and the powers of $Q$ can be related to each 
other by a pair of completely positive contractions $\theta: A\to B$ 
and $\theta_*: B\to A$ satisfying 
$\theta\circ\theta_*=\id _B$ (see Definition \ref{def2}).  
Significantly, when a $C^*$-dynamical system $(B,\alpha)$ is 
related in this way to $Q$ then it is 
uniquely determined by $Q$ up to conjugacy.  
We consider that a pair $(A,P)$ is 
{\em asymptotically stable} if there is a  
(necessarily unique) quasiautomorphism $Q: A\to A$ 
with the following property:
\begin{equation}\label{Eq03}
\lim_{n\to\infty}\|P^n(a)-Q^n(a)\|=0, \qquad a\in A.  
\end{equation}
See Theorem \ref{thm3}.  
Given the relation between $Q$ 
and $(B,\alpha)$ described above, 
one may conclude that when $(A,P)$ is 
stable in the sense of (\ref{Eq03}), 
then the {\em asymptotic} properties of the sequence 
$P, P^2, P^3,\dots$ are identical
with the asymptotic properties of the sequence 
of automorphisms $\alpha, \alpha^2,\alpha^3,\dots$.  
In particular, once one knows the $C^*$-dynamical 
system $(B,\alpha)$, one knows everything about the 
asymptotics of $(A,P)$.  This stability result for 
completely positive maps on \cstar s 
generalizes the classical Perron-Frobenius 
theorem on the structure of square matrices with 
nonnegative entries.  The connection between 
the stability assertion 
(\ref{Eq03}) and the Perron-Frobenius theorem 
is discussed more fully in Remark \ref{rem2}.  

The appropriate formulation of asymptotic 
stability for normal completely positive maps $P: M\to M$ on von Neumann 
algebras differs significantly from the formulation 
\cstar s, since it involves elements 
$\rho\in M_*$ of the predual 
rather than elements $a\in M$.  The appropriate 
formulation is this:  There should be a unique 
{\em normal} quasiautomorphism $Q: M\to M$ with the property 
that for every normal linear functional $\rho\in M_*$, one has 
\begin{equation}\label{Eq05}
\lim_{n\to\infty}\|\rho\circ P^n-\rho\circ Q^n\|=0.  
\end{equation}
In this case, there is a $W^*$-dynamical system 
$(N,\alpha)$ associated with $Q$ as in the case of 
\cstar s, and which is unique up to conjugacy.  

With this formulation of stability for 
von Neumann algebras, simple mixing 
of the type (\ref{Eq02}) becomes a special case 
of (\ref{Eq05}) as follows.  
Let $P:\mathcal B(H)\to \mathcal B(H)$ be a normal 
completely positive unit-preserving map for which 
there is a normal state $\omega$ on $\mathcal B(H)$ 
satisfying the strong mixing requirement 
(\ref{Eq02}).  Let $Q$ be the normal 
map of $\mathcal B(H)$ defined by $Q(x)=\omega(x)\mathbf 1$, 
$x\in \mathcal B(H)$.  Then $Q$ is a quasiautomorphism 
with range $\mathbb C\cdot\mathbf 1$, 
having the property that for every normal state 
$\rho$ of $\mathcal B(H)$, 
$$
\rho\circ Q^n=\omega, \qquad n=1,2,\dots.  
$$
It follows that (\ref{Eq02}) and (\ref{Eq05}) make 
the same assertion in this case.  
The $W^*$-dynamical 
system associated with this $Q$ is the trivial one 
$(\mathbb C, \id)$, $\id$ denoting the identity 
automorphism of the one-dimensional von Neumann 
algebra $\mathbb C$.  

These developments rest on 
some very general results for contractions 
acting on Banach spaces, and logic requires 
that we first work out 
this basic material  
in Section \ref{S:contr}.  We discuss the  
properties of quasiautomorphisms in 
Section \ref{S:paut}, and then give applications 
to \cstar s and von Neumann algebras 
in sections \ref{S:appl1}--\ref{S:appl2}, 
including examples.  
Similar results are valid for one-parameter semigroups, 
though beyond a few basic considerations in Section \ref{S:semig}, 
applications to semigroups are not developed 
here.  

The main hypothesis invoked 
in Theorems \ref{thm3} and \ref{thm4} is not necessary 
for the main conclusion concerning stability, and 
it is reasonable to ask if {\em every} 
normal unit-preserving completely positive 
map on $\mathcal B(H)$ is stable.  
The purpose of the last section is to show that there 
is a naturally-occurring class of such 
maps that are unstable.

\section{Locally Finite Contractions}\label{S:contr}

In this section we establish 
a general result in the 
category of Banach spaces, with contractions as maps.  
A {\em contraction} is an operator 
$T\in\mathcal B(X)$ acting on a complex Banach space $X$ 
that satisfies $\|T\|\leq 1$.  
We are concerned with the structure of contractions 
and with the 
asymptotic properties of their 
associated semigroups $\mathbf 1, T, T^2,\dots$.  
By an {\em automorphism} we mean an invertible isometry 
$U\in\mathcal B(X)$.  

\begin{rem}\label{rem1}
Given a pair of 
Banach spaces $X_1$, $X_2$ and an automorphism $U\in\mathcal B(X_1)$, 
there are many ways to introduce a 
norm on the algebraic direct sum $X_1\dotplus X_2$ 
so as to obtain a Banach space $X_1\oplus X_2$ with 
the property $\|x_1\|\leq \|x_1+x_2\|$ for all 
$x_k\in X_k$.  Settling on one of 
these norms, one can then define a 
contraction $T\in\mathcal B(X_1\oplus X_2)$ as 
the direct sum of operators 
$T=U\oplus \bf 0$, $\bf 0$ denoting  
the zero operator on $X_2$.  This is the 
most general example of a quasiautomorphism, 
a concept defined more concisely 
in operator-theoretic terms as follows.  
\end{rem}

\begin{defn}
A contraction $T\in\mathcal B(X)$ is called a 
{\em quasiautomorphism} if the restriction of $T$ 
to its range $TX$ is an 
automorphism of $TX$.  
\end{defn}

The range of a quasiautomorphism is necessarily closed, 
and every quasiautomorphism $T$ admits a unique 
``polar decomposition" 
$T=UE$, where $E$ is an idempotent contraction with range 
$TX$, and $U$ is an 
automorphism of $TX$.  
Indeed, $U$ is the restriction of $T$ to its 
range, and $E$ is the composition $U^{-1} T$.     
The projection $E$ commutes with 
$T$, and we have $T^n=U^nE$ for every $n=1,2,3,\dots$.  
Quasiautomorphisms are in an obvious sense 
the most rigid contractions.  

We emphasize that the term quasiautomorphism will be used 
below in other categories with more structure, 
and the attributes of quasiautomorphisms 
vary from one category to another.  
For example, when we work with 
completely positive maps 
on \cstar s in Section \ref{S:appl1}, 
quasiautomorphisms inherit  
the properties of maps in that category.  

The 
purpose of this section is to relate the {\em asymptotic} behavior of 
a broad class of contractions to that of quasiautomorphisms.   
These are the {\em locally finite} contractions, 
whose action on vectors is characterized 
as follows.  
The linear span of a finite set of vectors $x_1,\dots,x_n\in X$ 
is denoted by $[x_1,\dots,x_n]$.  

\begin{prop}\label{prop1}
Let $T$ be a contraction on a Banach space $X$.  For every 
vector $x\in X$, the following are equivalent:
\begin{enumerate}
\item[(i)]For every $\epsilon>0$ there is a positive 
integer $N$ such that 
$$
{\rm dist}(T^nx, [x, Tx, T^2x,\dots, T^Nx])\leq \epsilon, \qquad n=0,1,2,\dots.  
$$
\item[(ii)]
For every $\epsilon>0$ there is a finite-dimensional 
subspace $F\subseteq X$ such that 
$$
{\rm dist}(T^nx, F)\leq \epsilon, \qquad n=0,1,2,\dots.  
$$
\item[(iii)]The norm-closure of the orbit $\{x, Tx, T^2x,\dots\}$ is compact.  
\end{enumerate}
\end{prop}

\begin{proof}
The implication (i)$\implies$(ii) is trivial.    

(ii)$\implies$(iii): It suffices to show that for every $\epsilon>0$, the 
orbit of $x$ $O_x=\{x, Tx, T^2x,\dots\}$ can be covered 
by a finite union of balls of radius 
$\epsilon$.  To that end, fix $\epsilon>0$ and let $F$ be a finite dimensional 
subspace of $X$ such that every point of $O_x$ is within $\epsilon/2$ of 
$K=\{f\in F: \|f\|\leq \|x\|+1\}$.  Since $K$ is compact it may be covered 
by a finite union of balls $B_1\cup\cdots\cup B_r$ of radius at most $\epsilon/2$.  
After doubling the radius 
of each ball $B_i$, one obtains a finite union of balls of 
radius $\epsilon$ that covers $O_x$.  

(iii)$\implies$(i):  Fix $\epsilon>0$.  
Since $O_x$ is dense in its closure and its closure is
compact, there is an $N\geq 1$ such that every point of $O_x$ is 
within $\epsilon$ of $\{x, Tx, T^2x,\dots, T^Nx\}$, and \
(i) follows.    
\end{proof}

\begin{defn}\label{def3} A contraction $T$ acting on a 
Banach space $X$ is called {\em locally finite} 
if every vector  $x\in X$ satisfies the 
conditions of Proposition \ref{prop1}.  
\end{defn}

\begin{rem}
A vector $x\in X$ that is {\em algebraic} in the sense that 
$p(T)x=0$ for some nonzero polynomial $p(z)$ obviously satisfies 
condition (i).  
A vector $x$ will satisfy (ii) when it remains localized under 
the action of the nonnegative powers of $T$ in the sense that no 
subsequence of $x, Tx, T^2x,\dots$ can wander in an essential 
way through infinitely many dimensions.  

A straightforward argument shows that in general, 
the set of vectors 
that satisfy condition (iii) of 
Proposition \ref{prop1} 
is a closed linear 
subspace of $X$ that is invariant under the set of 
all operators in $\mathcal B(X)$ that commute 
with $T$.  It follows that $T$ will be locally finite if, for example, 
the set of all algebraic vectors has $X$ 
as its closed linear span.  
\end{rem}

We fix attention on eigenvectors $x$ of $T$ whose 
eigenvalues have maximum absolute value: $Tx=\lambda x$, where 
$|\lambda|=1$.  
Such an 
$x$ is called a {\em maximal eigenvector}.  
The {\em point spectrum} of an operator $T$ is the 
set of all eigenvalues of $T$, written $\sigma_p(T)$, 
and of course the point spectrum can be empty.  
An invertible 
isometry $U\in\mathcal B(X)$ is said to be 
{\em diagonalizable} if $X$ is spanned by the 
(necessarily maximal) eigenvectors of $U$, and 
for such operators $\sigma_p(U)$ is dense 
in the spectrum of $U$.  Our use of the 
term {\em diagonalizable} is not universal; 
for example, diagonalizable unitary 
operators are often said to have 
{\em pure point spectrum}.  
Nevertheless, this terminology will be convenient.  
The asymptotic behavior of 
locally finite contractions is described as follows:

\begin{thm}\label{thm1}
For every locally finite contraction $T$ acting on a Banach space $X$ there is a 
unique quasiautomorphism $S\in\mathcal B(X)$ such that 
\begin{equation}\label{Eq00}
\lim_{n\to\infty}\|T^nx-S^nx\|=0,\qquad x\in X.  
\end{equation}

The restriction $U$ of $S=UE$ to its range 
is diagonalizable, and we have 
\begin{equation}\label{Eq0}
\sigma_p(U)=\sigma_p(T)\cap \mathbb T,   
\end{equation}
where $\mathbb T$ is the unit circle.  
The projection $E$ is characterized as the unique 
idempotent in the set $\mathcal L$ of strong limit points 
of the powers of $T$
\begin{equation}\label{Eq000}
\mathcal L=\bigcap_{n=1}^\infty 
\left\{T^n, T^{n+1},T^{n+2},\dots\right\}^{-{\rm strong}}, 
\end{equation}
and $U$ is the restriction of $T$ to $EX$.  
\end{thm}

\begin{rem}
(\ref{Eq00}) asserts that a locally finite contraction has the same 
asymptotic behavior as an automorphism.  According 
to (\ref{Eq0}), the point spectrum of $U$ consists of all 
eigenvalues of $T$ that are associated with maximal eigenvectors.  
A consequence of the characterization (\ref{Eq000}) is that 
the projection $E$ will share the salient features of 
$\{T^n:n\geq 1\}$.  For example, when $T$ is a completely positive 
contraction acting on a \cstar\ whose powers do {\em not} tend 
to zero in the strong operator topology, 
then $E$ will be a completely 
positive idempotent of norm $1$.  
\end{rem}

The proof of Theorem \ref{thm1} will 
make use of the following known result 
from the theory of almost periodic 
representations of groups -- in our case the group is $\mathbb Z$.   
That material generalizes work of Harald Bohr (for the group $\mathbb R$  
\cite{BohrAP}) 
to arbitrary groups, the generalization being due 
to von Neumann and 
others (see pp. 245--261 of \cite{HewRoss1}, 
and pp. 310--312 of \cite{HewRoss2}).  

\begin{lem}\label{lem1}
Let $U$ be an invertible isometry acting on a Banach space $X$, 
and suppose that 
the $\mathbb Z$-orbit 
$\{U^nx: n\in\mathbb Z\}$
of every vector $x\in X$  is relatively norm-compact.  Then $U$ 
is diagonalizable.  
\end{lem}

We also require the following observation.  

\begin{lem}\label{lem2}
Let $T$ be a contraction on a Banach space $X$ such that 
the identity operator belongs to the strong closure of 
$\{T, T^2, T^3,\dots\}$.  Then $T$ is an 
automorphism of $X$.  
\end{lem}

\begin{proof}
For each $x\in X$ there is a sequence $n_1, n_2, \dots$ of positive 
integers such that $T^{n_k}x\to x$, and therefore 
$\|T^{n_k}x\|\to \|x\|$, as $k\to\infty$.  Since $n_k\geq 1$ for each $k$ it follows 
that $\|Tx\|\geq \|x\|$, hence $\|Tx\|=\|x\|$.  Similarly, 
$x$ belongs to the closure of $TX$.  These observations show that 
$T$ is an isometry with dense range, hence it is invertible.  
\end{proof}

\begin{proof}[Proof of Theorem \ref{thm1}]
Let $G\subseteq X$ be the (conceivably empty) 
set of all maximal eigenvectors of $T$ 
and let $\tilde T$ be the restriction of $T$ to $M=\overline{\rm span}\,G$.  

We claim first that $\tilde T$ is an invertible isometry.  According 
to Lemma \ref{lem2}, that will follow if 
we prove that the identity operator $\mathbf 1_M$ of $M$ belongs to the 
strong closure of $\{\tilde T, \tilde T^2,\tilde T^3,\dots\}$.  
For that, let 
$\{x_1,\dots, x_r\}$ be a finite subset of $G$.  It suffices 
to show that there is an increasing 
sequence $n_1<n_2<\dots$ of integers such that 
\begin{equation}\label{Eq2}
\lim_{k\to\infty}\|T^{n_k}x_i-x_i\|=0, \qquad i=1,\dots,r.  
\end{equation}
Noting that $Tx_i=\lambda_ix_i$ for $\lambda_1,\dots,\lambda_r\in\mathbb T$, 
we make use of a familiar 
result from Diophantine analysis which asserts that 
for every finite choice of elements $\lambda_1,\dots, \lambda_r$ in the 
multiplicative group $\mathbb T$, 
there is an increasing sequence $n_1<n_2<\dots\in\mathbb N$ such that  
$$
\lim_{k\to\infty}\lambda_i^{n_k}=1,\qquad i=1,\dots,r.  
$$
This sequence $n_1, n_2,\dots$ obviously satisfies (\ref{Eq2}).  
Thus Lemma \ref{lem2} implies that $\tilde T$ is an invertible 
isometry.  

Now let $N$ be the asymptotic null space 
$$
N=\{x\in X: \lim_{n\to\infty}\|T^nx\|=0\}.  
$$
We will show that $N$ and $M$ are complementary 
subspaces in the sense that $N\cap M=\{0\}$ and $N+M=X$.

Indeed, 
\begin{equation}\label{Eq3}
N\cap M=\{0\}  
\end{equation} 
follows from the preceding paragraph.  For if if $z$ is a vector in 
$N\cap M$ then since $T$ restricts to an isometry on $M$ we have 
$\|z\|=\|T^nz\|$ for every $n=1,2,\dots$, while since $z\in N$ 
we have $\|T^nz\|\to 0$ as $n\to \infty$.  Hence $z=0$.  

For every $x\in X$ consider the set of limit points 
$$
K_\infty(x)=\bigcap_{n=1}^\infty\overline{\left\{T^nx, T^{n+1}x, T^{n+2}x, \dots\right\}},   
$$
the bar denoting closure in the norm of $X$.  
We claim that $K_\infty(x)\subseteq M$ for every $x\in X$.  To prove that, 
fix $x$ and choose 
$z\in K_\infty(x)$.  We claim first that there is a sequence 
$n_1<n_2<\dots$such that 
\begin{equation}\label{Eq4}
\lim_{k\to \infty}T^{n_k}z=z, 
\end{equation}
the convergence being in norm.  
Indeed, by definition 
of $K_\infty(x)$ there is a sequence $m_1<m_2<\dots$ such that 
$T^{m_k}x$ converges to $z$.  We may assume that the $m_k$ increase 
as rapidly as desired by passing to a subsequence, and 
we choose $m_k$ so that the sequence of differences 
$n_k=m_{k+1}-m_{k}$ increases to $\infty$.  
Writing $z=T^{m_k}x+(z-T^{m_k}x)$ and estimating in the obvious 
way, we obtain
\begin{align*}
\|T^{n_k}z-z\|&\leq\|T^{n_k}T^{m_k}x-T^{m_k}x\| + 2\|z-T^{m_k}x\|\\
&=
\|T^{m_{k+1}}x-T^{m_k}x\| + 2\|z-T^{m_k}x\|, 
\end{align*}
hence $\|T^{n_k}z-z\|\to0$ as $k\to \infty$.  
Keeping $z\in K_\infty(x)$ fixed and 
choosing a sequence $n_1<n_2<\dots$ satisfying (\ref{Eq4}), we consider 
the set of vectors 
$$
M_z=\{y\in X: \lim_{k\to\infty}T^{n_k}y=y\}.  
$$
$M_z$ is a closed linear subspace of 
$X$ that contains $z$, and it is invariant 
under all operators that commute with $T$.  
We claim that the restriction 
$\tilde T$ of 
$T$ to $M_z$ is an invertible isometry.  Indeed, from 
the definition of $M_z$ it follows that $\tilde T^{n_k}$ 
converges strongly to the identity operator of $M$ as 
$k\to \infty$, hence the assertion follows after another
application of Lemma \ref{lem2}.  
Noting that $\tilde T^{-k}$ belongs to the strong closure 
of $\{\tilde T, \tilde T^2,\tilde T^3,\dots\}$ for every 
integer $k\geq 0$, it follows that $\tilde T$ satisfies 
the hypotheses of Lemma \ref{lem1}.  We conclude 
from that result 
that $M_z$ is spanned 
by maximal eigenvectors of $T$, and is therefore a subspace 
of $M$.  In particular, $z\in M$.

We show now that $N+M=X$.  Choose $x\in X$.  We will 
exhibit a vector $e\in M$ such that $\|T^nx - T^ne\|\to 0$ 
as $n\to \infty$.  
Indeed, the distance from $T^nx$ to $K_\infty(x)$ must 
decrease to zero as $n\to \infty$
because $K_\infty(x)$ is the intersection 
of the decreasing sequence of {\em compact} sets 
$\{T^nx, T^{n+1}x, T^{n+2}x,\dots\}^-$.  
Thus there is a sequence $k_n\in K_\infty(x)$ with the property 
$\|T^nx-k_n\|\to 0$ as $n\to \infty$.  We have shown above that 
$K_\infty(x)\subseteq M$ and that the restriction of $T$ to $M$ 
is an isometry.  It follows the restriction of $T$ to 
the compact metric space $K_\infty(x)$ defines an isometry of metric 
spaces.  Since by the definition of $K_\infty(x)$ it is clear 
that $TK_\infty(x)$ is dense in $K_\infty(x)$, it follows that 
$TK_\infty(x)=K_\infty(x)$.  So for every $n$ there 
is an element $\ell_n\in K_\infty(x)$ such that $k_n=T^n\ell_n$ 
and $\|k_n\|=\|\ell_n\|$, hence 
\begin{equation}\label{Eq5}
\lim_{n\to \infty}\|T^nx-T^n\ell_n\|=0.  
\end{equation}
Finally, by compactness of $K_\infty(x)$ there is a subsequence 
$m_1<m_2<\dots$ such that $\ell_{m_k}\to e\in K_\infty(x)$ as 
$k\to \infty$.  From (\ref{Eq5}) we deduce that 
\begin{equation*}
\lim_{k\to\infty}\|T^{m_k}x-T^{m_k}e\|=
\lim_{k\to\infty}\|T^{m_k}x-T^{m_k}\ell_{m_k}\|=0.
\end{equation*}
It follows that $\lim_{k\to\infty}\|T^{m_k}(x-e)\|=0$.  
This implies that $x-e\in N$ because the sequence of norms 
$\|T^n(x-e)\|$ decreases with $n$.  Thus 
$x=(x-e)+e$ is exhibited as an element of $N+M$.  

Let $E$ be the idempotent defined by 
$E\restriction_N=0$ and $Ey=y$ for $y\in M$.  
We claim that 
$E$ is a strong limit point of the sequence 
$\{T, T^2, T^3,\dots\}$ of powers of $T$.  For that, 
it suffices to show that for every integer $N\geq 1$, 
every $\epsilon>0$, every finite set $x_1,\dots, x_n$ of 
maximal eigenvectors of $T$, and every finite set 
$z_1,\dots, z_m\in N$, there is an integer $p\geq N$ 
such that 
\begin{equation}\label{Eq6}
\|T^px_k-x_k\|\leq \epsilon \quad {\rm and }\quad  \|T^pz_j\|\leq \epsilon,
\qquad 1\leq k\leq n,\quad 1\leq j\leq m.  
\end{equation}
Writing $Tx_k=\lambda_kx_k$, $1\leq k\leq n$, the Diophantine approximation 
employed above shows that there is an infinite set of positive 
integers $p$ such that 
$|\lambda_1^p-1|\leq \epsilon, \dots, |\lambda_n^p-1|\leq \epsilon$.  
Since $\|T^pz_1\|,\dots,\|T^pz_m\|$ all tend to zero with large 
$p$, it is apparent that we can satisfy (\ref{Eq6}) with infinitely 
many values of $p$.  This shows that $E$ is a strong cluster 
point of $\{T^n:n\geq1\}$, and in particular  $\|E\|\leq 1$.  

If we set $S=TE$, then we may conclude 
from the preceding discussion 
that $S$ is a quasiautomorphism satisfying 
both (\ref{Eq00}) and (\ref{Eq0}).

We claim now that $E$ is the only idempotent that 
can be a strong cluster point of $\{T^n: n\geq 1\}$.  Let 
$F$ be such another such limit point.  Since both $E$ and $F$ 
are idempotents, to show that $F=E$ 
it suffices to show that $\ker F\subseteq\ker E=N$ and $FX\subseteq EX=M$.  
If z is any vector in the kernel of $F$, then there is a sequence 
$n_k\to\infty$ such that $T^{n_k}z\to Fz=0$, hence 
$\|T^{n_k}z\|\to 0$, as $k\to \infty$.  The latter implies 
that $\lim_{n\to \infty}\|T^nz\|=0$ since the norms 
$\|T^nz\|$ decrease with $n$, hence $z\in N$.  If $y$ is a 
vector in the range of $F$ then there is a sequence 
$m_k\to\infty$ such that $T^{m_k}y\to Fy=y$, and 
in particular $y\in K_\infty(y)$.  
We have already proved that $K_\infty(y)\subseteq M$, 
and therefore $y\in M=EX$.   

It remains to show that $S=TE$ is the only quasiautomorphism 
in $\mathcal B(X)$ that satisfies (\ref{Eq00}).
Let $R=UF$ be the polar decomposition of 
another quasiautomorphism such 
that $\lim_{n\to\infty}\|S^nx-R^nx\|=0$ 
for every $x\in X$.  We claim that $F=E$ and $U=T\restriction_M$.    
Since both $E$ and $F$ are 
idempotents, the first assertion will follow if we 
show that $N=\ker E\subseteq \ker F$ and that 
$M=EX\subseteq FX$.  Indeed, if $x\in \ker E$ then 
(\ref{Eq00}) implies that 
$$
\|Fx\|=\lim_{n\to\infty}\|U^nFx\|=\lim_{n\to\infty}\|T^nEx-U^nFx\|=
\lim_{n\to\infty}\|S^nx-R^nx\|=0,
$$
hence $x\in\ker F$.  $M\subseteq FX$ will follow if we 
show that every maximal eigenvector $x$ for $T$ belongs 
to the range of $F$.  Writing $Tx=Sx=\lambda x$ for some $\lambda\in\mathbb T$, 
we have $\|x-\bar\lambda^nR^nx\|=\|S^nx-R^nx\|\to0$ as $n\to\infty$, 
which implies that the distance from $x$ to $FX=RX$ is zero.  
Finally, to show that $U=T\restriction_M$, choose $x\in FX=EX$ 
and write 
\begin{align*}
\|Tx-Ux\|&=\|U^n(Tx-Ux)\|=\|U^nTx-S^{n+1}x + (S^{n+1}x-U^{n+1}x)\|\\
&\leq 
\|R^nTx-S^nTx\| + \|S^{n+1}x-R^{n+1}x\|.  
\end{align*}
As $n\to\infty$, both terms on the right tend to zero by hypothesis.  
It follows that $\|Tx-Ux\|=0$, and therefore $U=T\restriction_{EX}$.  
\end{proof}

\begin{rem}
Perhaps it is worth pointing out that the set $\mathcal L$ of 
strong cluster points (\ref{Eq000}) is compact in its relative 
strong operator topology.  Indeed, $\mathcal L$ is a compact 
topological group with respect to operator 
multiplication, whose unit is  $E$.  Since 
we do not require this fact, we omit the proof.  
\end{rem}

\section{Quasiautomorphisms of \cstar s}\label{S:paut}

The notion of quasiautomorphism
must be interpreted appropriately when it is applied to completely 
positive maps on \cstar s.  The purpose of this section is to 
make some observations that show how a quasiautomorphism of a \cstar \ can be related to to an ordinary  
$*$-automorphism of a different \cstar;  
and that in fact the $C^*$-dynamical 
system associated with the 
quasiautomorphism is unique
up to conjugacy.

By a {\em CP contraction} we 
mean a completely positive 
linear map $P: A\to A$ defined on a \cstar\ $A$ such 
that $\|P\|\leq 1$.  If $A$ has a unit $\bf 1$, then 
a completely positive map $P: A\to A$ is a CP contraction iff 
$\|P({\bf 1})\|\leq 1$; but in general, we may speak 
of CP contractions even when $A$ fails to posses a 
unit.  
Throughout the section, $A$ will denote a \cstar.  We 
first show how, starting with an automorphism of another 
\cstar\ that is suitably related to $A$, one obtains 
a CP contraction on $A$ with special features.  
We use the term {\em $C^*$-dynamical system} 
to denote a pair $(B,\beta)$ consisting 
of a \cstar\ $B$ and a $*$-automorphism $\beta: B\to B$.  
Two $C^*$-dynamical systems $(B_1,\beta_1)$ and $(B_2,\beta_2)$ 
are said to be {\em conjugate} if there is a $*$-isomorphism 
$\theta: B_1\to B_2$ satisfying 
$\theta\circ\beta_1=\beta_2\circ\theta$.  

\begin{prop}\label{prop5}
Let $(B,\beta)$ be a $C^*$-dynamical system and 
let $\theta: A\to B$,  $\theta_*: B\to A$  be 
a pair of completely positive contractions 
satisfying $$\theta\circ\theta_*=\id _B.$$  
Let $P: A\to A$ be 
the CP contraction defined 
by 
\begin{equation}\label{Eq9}
P=\theta_*\circ\beta\circ\theta.
\end{equation}

Then $E=\theta_*\circ\theta$ is an idempotent 
CP contraction on $A$ with range $P(A)$,  
$\ker E=\ker P$, $P E=E P=P$, 
and the restriction $\alpha$ of 
$P$ to $E(A)$ is a surjective completely isometric map 
with the property $P^n=\alpha^n\circ E$, 
for every $n=1,2,\dots$.  

If $(\tilde B,\tilde\beta)$ is another $C^*$-dynamical system 
that is similarly related to $P$, 
$P=\tilde\theta\circ\tilde \beta\circ\tilde\theta_*$, 
where $\tilde\theta: A\to\tilde B$ and $\tilde\theta_*: \tilde B\to A$ 
are completely positive contractions with 
$\tilde\theta\circ\tilde\theta_*=\id_{\tilde B}$, then 
the $C^*$-dynamical systems 
$(\tilde B,\tilde\beta)$ and $(B,\beta)$ are 
naturally conjugate.  
\end{prop}

\begin{proof}
Since $\theta\circ\theta_*=\id _B$, it follows that $\theta$ 
is surjective, $\theta_*$ is injective, and 
$$
E^2=\theta_*\circ\theta\circ\theta_*\circ\theta=
\theta_*\circ\theta=E.     
$$  
Since $\theta$ and $\beta\circ\theta$ are both surjective, 
\begin{equation}\label{Eq1}
\ran P=\theta_*(\beta(\theta(A)))=\theta_*(B)=\theta_*(\theta(A))=\ran E;    
\end{equation}
and since $\theta_*$ and $\theta_*\circ\beta$ are 
both injective, 
\begin{equation}\label{Eq1a}
\ker P=\ker \theta_*\circ\beta\circ\theta=
\ker\theta=\ker\theta_*\circ\theta=\ker E.  
\end{equation}

Similarly, one verifies directly that 
the restriction $\alpha$ of $P$ to $E(A)$ satisfies 
$\alpha^n(E(a))=\theta_*\circ\beta^n\circ\theta(E(a))$ for 
$a\in A$, $n=1,2,\dots$.  In particular, 
$\alpha$ is the restriction of $\theta_*\circ\beta\circ\theta$ 
to $E(A)$, a completely isometric surjective map 
of $E(A)$ onto itself.  

Suppose that $(\tilde B, \tilde \beta)$ is another $C^*$-dynamical 
system and $\tilde\theta: A\to B$ and $\tilde\theta_*: B\to A$ 
are CP contractions satisfying 
$\tilde\theta\circ\tilde\theta_*=\id _{\tilde B}$ and 
$P=\tilde\theta_*\circ\tilde\beta\circ\tilde \theta$.  In this 
case we have a second CP idempotent $\tilde E: A\to A$ defined 
by $\tilde E=\tilde\theta_*\circ\tilde\theta$, and 
we claim that $\tilde E=E$.  Indeed, since both $\tilde E$ and 
$E$ are idempotents it suffices to show that they have 
the same kernel and the same range; and  
(\ref{Eq1}) and (\ref{Eq1a}) imply that 
$\ker\tilde E=\ker P=\ker E$ and  
$\ran\tilde E=\ran P=\ran E$.

We define CP contractions $\phi: B\to \tilde B$ 
and $\tilde\phi: \tilde B\to B$ by 
$$
\phi = \tilde\theta\circ\theta_*,\qquad \tilde\phi
=\theta\circ\tilde\theta_*.  
$$
We have 
$$
\phi\circ\tilde\phi=
\tilde\theta\circ \theta_*\circ\theta\circ\tilde\theta_*
=\tilde\theta\circ E\circ\tilde\theta_*=
\tilde\theta\circ\tilde E\circ\tilde\theta_*=
\id_{\tilde B}^2=\id_{\tilde B}.  
$$
Similarly, $\tilde\phi\circ\phi=\id_B$, so that the 
maps $\phi$, $\tilde\phi$ are completely isometric 
completely positive 
maps that are inverse to each other.  
Since $B$ and $\tilde B$ are both \cstar s, 
we may conclude that 
$\phi$ is a $*$-isomorphism of $B$ onto $\tilde B$ 
with inverse $\tilde\phi$.  

It remains to show that $\phi\circ\beta=\tilde\beta\circ\phi$.  
Composing the identity 
$$
\theta_*\circ\beta\circ\theta=
\tilde\theta_*\circ\tilde\beta\circ\tilde\theta=P
$$ 
on the left with $\tilde\theta$ gives 
$$
\phi\circ\beta\circ\theta=
\tilde\theta\circ\tilde\theta_*\circ\tilde\beta\circ\tilde\theta
=\tilde\beta\circ\tilde\theta, 
$$
and after composing with $\theta_*$ on the right 
we obtain $\theta\circ\beta=\tilde\beta\circ\theta$.  
\end{proof}

\begin{prop}\label{prop6}
For every CP contraction $P: A\to A$ on a \cstar\ 
$A$, the following are equivalent.  
\begin{enumerate} 
\item[(i)]$P$ admits a factorization $P=\alpha\circ E$ 
where $E: A\to A$ is an idempotent completely positive contraction and 
$\alpha$ is a completely isometric linear map of $E(A)$ 
onto itself.  
\item[(ii)]There is a $C^*$-dynamical system $(B,\beta)$ 
that is related to $P$ as in (\ref{Eq9}).  
\end{enumerate}

If $A$ has a unit $\bf 1$ and $P(\bf 1)=\bf 1$, then 
(i) can be replaced with 
\begin{enumerate}
\item[(i)$^\prime$] The restriction of $P$ to $P(A)$ is 
a surjective complete isometry.  
\end{enumerate}
\end{prop}

\begin{proof}The implication (ii)$\implies$(i) follows 
from Proposition \ref{prop5}.  

(i)$\implies$(ii): The hypothesis (i) obviously 
implies that $EP=P=PE$.  

A result of Choi and Effros 
\cite{ChE} implies that $E(A)$ is a \cstar\ 
with respect to the multiplication defined on it by 
$x\bullet y=E(xy)$, for $x,y\in E(A)$ (one uses the norm of $A$ 
and the vector space operations and $*$-operation inherited 
from $A$).  Let $B$ be this \cstar.  

We may consider 
$E$ as a completely positive contraction of $A$ onto 
$B$; let $\theta$ be that map, and let $\theta_*$ be 
the natural inclusion of $B=E(A)\subseteq A$.  
Obviously, $\theta\circ\theta_*=\id_B$.  
By hypothesis, the restriction of
$P$ to $E(A)$ is a surjective completely isometric map 
of the operator space $E(A)$ onto itself, 
and therefore it defines a completely isometric linear 
map $B$ onto itself, which we denote by $\beta$.  
We have to show that $\beta$ is a $*$-automorphism of 
$B$ and the maps $\theta$, $\theta_*$ relate $\beta$ 
to $P$ as in (\ref{Eq9}).  

We claim first that 
$\beta$ is also a positive linear map on $B$, 
i.e., $\beta(x^*\bullet x)\geq 0$ for every $x\in E(A)$.  
To see that, fix $x\in E(A)$, choose a positive 
linear functional $\rho$ on $B$, and consider 
the linear functional defined on $A$ 
by $\omega(a)=\rho(\theta(a))=\rho(E(a))$, $a\in A$.  
$\omega$ is a positive linear functional on 
$A$, hence 
$$
\rho(\beta(x^*\bullet x)=\rho(P(E(x^*x)))=
\rho(E(P(x^*x)))=\omega(P(x^*x))\geq0.  
$$
Since $\rho$ is an arbitrary positive linear functional 
on $B$,  $\beta(x^*\bullet x)\geq 0$ follows, 
hence $\beta$ is a positive linear map.  An obvious 
variation of this argument (that we omit) shows that $\beta$ induces 
a positive linear map on every matrix algebra 
$M_n\otimes B$ over $B$, hence $\beta$ is a completely positive 
linear map that is also completely isometric.  It follows 
that $\beta$ is a $*$-automorphism of $B$.  

Finally, to check that (\ref{Eq9}) is satisfied, we have 
$$
\theta_*\circ\beta\circ\theta(a)=P(E(a))=P(a),\qquad a\in A, 
$$
since $P\circ E=\alpha\circ E=P$.  

Finally, assuming that $A$ has a unit $\mathbf 1$ and 
$P(\mathbf 1)=\mathbf 1$, then $E(A)$ is an operator 
system, and a unit-preserving linear map of 
one operator 
system to another 
is completely positive iff it is completely 
contractive \cite{arvSubalgI}.  So in this case 
(i)$^\prime$ is equivalent to (i).  
\end{proof}

\begin{defn}\label{def2}
A completely positive contraction $P: A\to A$ is called 
a {\em quasiautomorphism} if the conditions of Proposition 
\ref{prop6} are satisfied.  
\end{defn}

Some concrete examples of quasiautomorphisms 
are given in Section \ref{S:appl1}.  
The preceding remarks support the following 
point of view: {\em the nontrivial 
behavior of the powers of a quasiautomorphism is identical with the behavior 
of the powers of a uniquely determined 
automorphism of a \cstar}.

\section{Applications to \cstar s}\label{S:appl1}

In this section we 
describe an application of Theorem \ref{thm1} 
to completely positive maps on 
\cstar s and describe how that result 
provides a noncommutative generalization of 
the Perron-Frobenius theorem.  
We also exhibit a variety of examples 
of locally finite completely positive maps 
on infinite-dimensional \cstar s.  

The term {\em locally finite} applies to 
completely positive contractions  
exactly as stated in Definition \ref{def3}, 
and of course in this context the strong operator topology 
is the topology of point-norm convergence: a net of linear 
maps $L_i: A\to A$ converges strongly to a linear map $L: A\to A$ iff 
one has 
$$
\lim_{i\to\infty}\|L_i(a)-L(a)\|=0, \qquad a\in A.  
$$

\begin{thm}\label{thm3}
Let $A$ be a \cstar\ and let $P: A\to A$ be 
a locally finite 
completely positive contraction.  
Then there is a unique quasiautomorphism $Q=\alpha\circ E$ of $A$ such that 
\begin{equation}\label{Eq13}
\lim_{n\to\infty}\|P^n(a)-Q^n(a)\|=0, \qquad a\in A.  
\end{equation}

The completely positive idempotent 
$E$ is characterized as the unique 
idempotent in the set of strong cluster points of 
$\{P, P^2, P^3,\dots\}$, $\alpha$ is 
the restriction of $P$ to the operator space $E(A)$, 
and $E(A)$ is the norm-closed linear span of the set 
of all maximal eigenvectors of $P$.  
\end{thm}

\begin{proof}
By Theorem \ref{thm1}, there is a unique idempotent 
$E$ in the set of strong limit points of 
$\{P, P^2,\dots\}$.  
Being a limit in the strong operator 
topology of a net of powers of 
$P$, $E$ is a completely positive contraction.  
Moreover, the general results of 
Theorem \ref{thm1} imply that 
$P$ restricts to an isometry $\alpha$ of $E(A)$ onto 
itself in such a 
way that the map $Q=\alpha\circ E$ satisfies 
$$
\lim_{n\to\infty}\|P^n(a)-Q^n(a)\|=0, \qquad a\in A.  
$$

We claim that $\alpha$ is a {\em completely} isometric 
linear map of $E(A)$ onto itself.  Indeed, for 
fixed $n=2,3,\dots$, consider the map of $M_n\otimes A$ 
defined on $n\times n$ matrices over $A$ by 
$$
\id_n\otimes P:(a_{ij})\mapsto (P(a_{ij})).  
$$
$\id_n\otimes P$ satisfies the same hypotheses as 
$P$, and since $\id_n\otimes \alpha$ 
is the restriction of $\id_n\otimes P$ to $M_n\otimes E(A)$, 
we may argue exactly as above to conclude 
that $\id_n\otimes \alpha$ is isometric 
on $\id_n\otimes E(A)$.  Hence $\alpha$ is completely 
isometric.  

It follows that $Q=\alpha\circ E$ is a quasiautomorphism 
of $A$ in the sense of Definition \ref{def2}.  The 
remaining assertions, including the uniqueness of 
$Q$, now follow from Theorem \ref{thm1}.  
\end{proof}

\begin{rem}[Relation to the Perron-Frobenius theory]
\label{rem2}
Frobenius' generalization \cite{frob} of 
Perron's theorem \cite{perron1}\cite{perron2} 
on square 
matrices with positive entries can be viewed 
as a result that provides information about the structure 
and properties of positive linear maps acting on 
finite-dimensional commutative \cstar s.  Indeed, 
every $n\times n$ matrix with nonnegative entries 
acts naturally on complex column vectors 
as a positive linear 
map, and every positive linear map of $\mathbb C^n$ 
arises in that way.  Recall too that a positive 
linear map on a commutative \cstar\ is 
automatically completely positive.  

In order to 
simplify the following remarks, we start with a positive 
linear map $P: A\to A$ on a finite-dimensional 
commutative 
\cstar\ $A$ satisfying $P(\bf 1)=\bf 1$, 
in which case both the norm and 
spectral radius of $P$ are $1$.   Thus the 
first assertion of the Perron-Frobenius theorem, 
namely that there is a nonzero positive element 
of $A$ that is fixed under $P$, is automatic.  For 
purposes of this discussion, the principal assertions of 
Frobenius' result (page 65 of \cite{gant}) can be paraphrased 
for positive maps as follows.  

\newtheorem*{nonum}{Perron-Frobenius Theorem}
\begin{nonum}
Assume further that $P$ is 
{\em irreducible} in the sense that the 
only projections $e\in\mathcal A$
satisfying $P(e)\leq e$ are $e=\bf 0$ and $e=\bf 1$, 
and let $\{\lambda_0,\dots,\lambda_{k-1}\}$ be the 
distinct eigenvalues 
of $P$ that lie on the unit circle, $1\leq k\leq \dim A$.  

Then each $\lambda_j$ is a simple eigenvalue and 
$\lambda_0,\dots, \lambda_{k-1}$ 
are the distinct $k$th roots of unity; hence we can 
arrange that $\lambda_j=\zeta^j$, where $\zeta=e^{2\pi i/k}$.  
Moreover, the matrix of $P$ relative to a 
basis of minimal projections of $A$ has the form $UCU^{-1}$, 
where $U$ is a permutation matrix 
and $C$ is a ``cyclic" 
matrix of rectangular blocks
\begin{equation*}
C=
\begin{pmatrix}
0&C_0&0&\dots&0\\
0&0&C_1&\dots&0\\
\hdotsfor[2.0]{5}\\
0&0&0&\dots&C_{k-2}\\
C_{k-1}&0&0&\dots&0
\end{pmatrix}, 
\end{equation*}
in which the diagonal blocks are square.  

\end{nonum}

We now describe how these 
combined assertions about the structure of $P$ fit 
naturally into the context of Theorem \ref{thm3}.  The 
displayed cyclic structure of $C$, together with the 
fact that $P=UCU^{-1}$, implies that there is a set of 
mutually orthogonal projections $e_0,\dots,e_{k-1}\in A$, 
with $e_0+\cdots+e_{k-1}=\bf 1$, which are permuted 
cyclically by $P$ in the sense that 
\begin{equation}\label{Eq17}
P(Ae_i)\subseteq Ae_{i\dot+1}, 
\end{equation}
where $\dot+$ denotes addition modulo $k$.  Let 
$B$ be the $C^*$-subalgebra of $A$ spanned by the 
projections $e_0,\dots, e_{k-1}$.  
(\ref{Eq17}) implies that $P(e_i)\leq e_{i\dot+1}$, 
and after summing on $i$ we find that 
equality must hold for each $i$ 
because  $e_0+\cdots+e_{k-1}=\bf 1$ 
and $P(\bf 1)=\bf 1$.  
Thus the restriction of $P$ to $B$ is the 
$*$-automorphism $\alpha$ of $B$ determined by 
$\alpha(e_i)=e_{i\dot+ 1}$, $0\leq i\leq k-1$.  

This automorphism $\alpha$ 
is the ``isometric" part of the quasiautomorphism $Q$ 
that is associated 
with $P$ by 
Theorem \ref{thm3}.   To see that, one first observes 
that $B$ is spanned by the set of maximal eigenvectors 
of $P$.  Indeed, an elementary 
argument shows that $B$ is spanned by the 
set of elements 
$$
x_\ell=e_0+\bar\zeta^\ell e_1+\bar\zeta^{2\ell}e_2+
\cdots+\bar\zeta^{(k-1)\ell}e_{k-1}, 
\qquad 0\leq \ell\leq k-1, 
$$
and clearly $P(x_\ell)=\zeta^\ell x_\ell$ for all $\ell$.  Since 
the eigenvalues
$\lambda_\ell=\zeta^\ell$ are all simple, it follows that 
$B=[e_0,\dots, e_{k-1}]=[x_0,\dots,x_{k-1}]$ is the space 
spanned by all maximal eigenvectors.  Thus, 
Theorem \ref{thm3} implies that 
the $*$-automorphism $\alpha=P\restriction_B$ is related to 
$Q$ by $Q=\alpha\circ E$ where $E$ is the unique idempotent 
limit point of $\{P^n\}$.  
We now identify $E$.  Since $B$ is spanned by 
the maximal eigenvectors of $P$, the proof of Theorem \ref{thm1} 
implies that we have a direct sum decomposition of 
finite-dimensional vector spaces 
$$
A=B\oplus\{z\in A: \lim_{n\to\infty}\|P^n(z)\|=0\}.    
$$
Thus the sequence $P^k, P^{2k}, P^{3k}, \dots$ converges 
to an idempotent with range $B$, and another application 
of Theorem \ref{thm3} shows that $E=\lim_n P^{nk}$.

In fact, given the relation between quasiautomorphisms and 
$*$ automorphisms in Proposition \ref{prop6}, 
it is not hard to turn this argument around to deduce 
Theorem 2 of \cite{gant} (including 
the permutation formula (\ref{Eq17}))
from Theorem \ref{thm3} above; 
and in this sense {\em one can regard Theorem \ref{thm3} as 
a generalization of the Perron-Frobenius theorem 
to \cstar s}.  
\end{rem}

Our search for a result like Theorem \ref{thm3} 
was inspired in part by a recent observation 
of Greg Kuperberg on the existence of idempotent 
limits of powers of a 
completely positive map on 
a finite-dimensional \cstar:   

\begin{thm}[Kuperberg]\label{kupThm}
Let $A$ be a finite-dimensional \cstar\ and let 
$P: A\to A$ be a unital completely positive map.  
There is a sequence of integers $0<n_1<n_2<\cdots$ 
such that $P^{n_k}$ converges to a unique completely 
positive idempotent map $E: A\to A$.  
\end{thm}
The uniqueness assertion 
means that $E$ does not depend on the sequence $n_k$ 
in the sense that if $m_1<m_2<\cdots$ is another increasing sequence 
for which the powers $P^{m_k}$ converge to an {\it idempotent} $F$, 
then $F=E$.  An elementary 
proof of Theorem \ref{kupThm} is sketched in \cite{kup1}.  

Of course, any unital 
completely positive linear map that acts on a finite dimensional 
\cstar\ must be locally finite, and therefore satisfies 
the hypotheses of Theorem \ref{thm3}.  But there are many 
others as well, and we 
now briefly describe some examples that act on 
familiar \cstar s.  

\vskip0.1in

{\bf Example.}  Let $G$ be a discrete group, let 
$U: G\to\mathcal B(\ell^2(G))$ be the regular representation of 
$G$ on its natural Hilbert space, 
and consider the reduced group \cstar\ $\mathcal A=C^*\{U_x: x\in G\}$ 
of $G$.  For every state $\rho$ of $\mathcal A$ there 
is a naturally associated positive definite function 
$\phi: G\to\mathbb C$, defined 
by $\phi(x)= \rho(U_x)$, and one has 
$|\phi(x)|\leq \phi(e)=\rho(\bf 1)=1$ for all $x\in G$.  
Notice first that the ``kernel" of $\phi$ 
$$
K=\{x\in G: |\phi(x)|=1\} 
$$
is a subgroup of $G$ and 
the restriction of $\phi$ to $K$ 
is a character of $K$.  Indeed, the GNS construction 
provides us with a unitary representation 
$V: G\to H$ and a cyclic vector $\xi$ for $V$ so that 
$$
\phi(x)=\langle V_x\xi,\xi\rangle, \qquad x\in G.  
$$
Noting that $\|V_x\xi-\phi(x)\xi\|^2=2 - 2|\phi(x)|^2$ for 
any $x\in G$, 
it follows that $|\phi(x)|=1$ iff $\xi$ is an eigenvector 
for $V_x$ in the sense that $V_x\xi=\phi(x)\xi$.   
Thus $K$ is a subgroup
on which $\phi$ is multiplicative.  

\begin{prop}\label{prop2}
For every state $\rho$ of $\mathcal A$, there is a unique 
completely positive linear 
map $P: \mathcal A\to \mathcal A$ satisfying 
\begin{equation}\label{Eq10}
P(U_x)=\rho(U_x)U_x,\qquad x\in G.     
\end{equation}
$P$ is a locally finite 
map with the following properties.  

Let $\mathcal B$ be the $C^*$-subalgebra of $\mathcal A$ 
generated by $\{U_x: |\rho(U_x)|=1\}$.  Then 
\begin{enumerate}
\item[(i)]
The restriction $\alpha$ 
of $P$ to $\mathcal B$ is a $*$-automorphism of $\mathcal B$.  
\item[(ii)]
There is a unique completely positive map $E$ defined on 
$\mathcal A$ by 
$$
E(U_x)=
\begin{cases}
U_x, &|\rho(U_x)|=1\\
0, & |\rho(U_x)|<1,    
\end{cases}
$$
and $E$ is an idempotent with range $\mathcal B$.  
\item[(iii)]
For every $A\in\mathcal A$ we have 
$$
\lim_{n\to\infty}\|P^n(A)-\alpha^n\circ E(A)\|=0.  
$$
\end{enumerate}
\end{prop}

\begin{proof}[Sketch of proof]To see that there 
is a completely positive map $P:\mathcal A\to \mathcal A$ 
satisfying (\ref{Eq10}), consider the unitary representation 
$$
W_x=U_x\otimes U_x \in\mathcal B(\ell^2(G)\otimes\ell^2(G)), 
\qquad x\in G.
$$ 
By Proposition 4.2 
of \cite{vallinHopf}, $W$ is weakly contained 
in the regular representation, so there is a representation 
$\pi: \mathcal A\to\mathcal B(\ell^2(G)\otimes \ell^2(G))$ 
satisfying $\pi(U_x)=U_x\otimes U_x$, $x\in G$.  Letting 
$Q:\mathcal A\otimes \mathcal A\to \mathcal A$ be the 
slice map $Q(A\otimes B)=\rho(A)B$, one finds that the composition 
$Q\circ\pi$ satisfies (\ref{Eq10}).  The proofs 
of the remaining assertions 
are straightforward.  
\end{proof}

I want to thank Marc Rieffel for the reference \cite{vallinHopf}.  
There are many variations of this example, including 
some natural examples acting on irrational rotation 
\cstar s.  
These examples all have the feature that the range 
$E(\mathcal A)$ of 
the completely positive idempotent 
$$
E\in\bigcap_{n\geq0}\overline{\{P^n, P^{n+1},\dots\}}^{\rm strong}
$$ 
is already $C^*$-subalgebra of $\mathcal A$.  
However, that is an 
artifact of this class of examples, and perhaps it is 
worth pointing out that in general, $E(\mathcal A)$ 
need not be a subalgebra of $\mathcal A$.  The following 
examples illustrate the point.

\vskip0.1in

{\bf Example.}
Let $\mathcal T$ be the Toeplitz \cstar, the 
\cstar \ generated by the simple unilateral shift.  
The familiar exact sequence of \cstar s
\begin{equation}\label{Eq11}
\begin{CD}
0 @>>> \mathcal K@>>> \mathcal T
@>>\pi> C(\mathbb T)@>>> 0
\end{CD}
\end{equation}
has a positive linear lifting $\pi_*: C(\mathbb T)\to \mathcal T$ 
which carries a symbol 
$f\in C(\mathbb T)$ to its associated Toeplitz operator 
$T_f$.  Every homeomorphism $h:\mathbb T\to \mathbb T$ gives 
rise to a $*$-automorphism 
$\beta$ of $C(\mathbb T)$ via $\beta(f)=f\circ h$, 
and after fixing $h$ one obtains a quasiautomorphism $Q: \mathcal T\to \mathcal T$ 
as in Proposition \ref{prop5},
$$
Q(A)=\pi_*(\beta(\pi(A))), \qquad A\in\mathcal T.  
$$
In more explicit terms, every operator in $\mathcal T$ 
admits a unique decomposition $A=T_f+K$, where $f\in C(\mathbb T)$, 
$K\in \mathcal K$ \cite{arvSpecTh}, and we have 
$$
Q(T_f+K)=T_{f\circ h}, \qquad f\in C(\mathbb T), \quad K\in \mathcal K.  
$$
This quasiautomorphism has ``polar decomposition" 
$Q=\alpha\circ E$, where 
$E: \mathcal T\to \mathcal T$ is the completely 
positive idempotent 
$$
E(T_f+K)=T_f,\qquad f\in C(\mathbb T),\quad K\in \mathcal K,   
$$
and $\alpha$ is the completely isometric linear map 
defined on $E(\mathcal T)$ by 
$$
\alpha(T_f)=T_{f\circ h},\qquad f\in C(\mathbb T).  
$$
Note that in these examples, 
$E(\mathcal T)$ is the space of all Toeplitz 
operators with continuous symbol, an operator system that is 
not a $C^*$-subalgebra of $\mathcal T$.  

As a variation on this example, 
let $P_0$ be a 
unital completely positive map that acts on a 
finite-dimensional \cstar\ 
$\mathcal A$, and let $P=P_0\otimes Q$ be 
the completely 
positive linear map on $\mathcal A\otimes \mathcal T$ 
that satisfies 
\begin{equation}\label{Eq14}
P(X\otimes A)=P_0(X)\otimes Q(A),\qquad X\in \mathcal A, 
\quad A\in\mathcal T.
\end{equation}
Since $P_0$ acts on a finite-dimensional \cstar\ it is 
locally finite, hence there is a unique quasiautomorphism 
$Q_0:\mathcal A\to \mathcal A$ satisfying 
$$
\lim_{n\to\infty}\|P_0^n(X)-Q_0^n(X)\|=0, \qquad 
X\in \mathcal A.  
$$  
It follows that 
$$
\lim_{n\to\infty}\|P^n(X\otimes A)-Q_0^n(X)\otimes Q^n(A)\|=0, 
\qquad X\in \mathcal A,\quad A\in\mathcal T.   
$$
Thus, $Q_0\otimes Q$ is 
the quasiautomorphism of $\mathcal A\otimes \mathcal T$ 
that is asymptotically associated 
with $P$.  The range of $Q_0\otimes Q$ is 
certainly an operator system, but 
it is never a $C^*$-subalgebra of 
$\mathcal A\otimes \mathcal T$.  

It goes without saying 
that one can obtain a great variety of such examples 
by replacing the Toeplitz diagram 
(\ref{Eq11}) with other linearly split short exact sequences 
of \cstar s.  

\section{Applications to von Neumann algebras}\label{S:appl2}

We now describe how Theorem \ref{thm3} must be modified 
for normal completely positive maps acting 
on von Neumann algebras.

\begin{defn}
By a {\em quasiautomorphism} of a von Neumann algebra $M$ we mean 
a {\em normal} unit-preserving completely positive linear map $P:M\to M$, 
such that $P(M)$ is norm-closed and 
$P$ restricts to a completely 
isometric linear map of $P(M)$ onto itself.  
\end{defn}

\begin{rem}[Structure of Quasiautomorphisms]\label{rem3}
Note first that the range $P(M)$ 
of a quasiautomorphism must be weak$^*$-closed.  
Indeed, since $P$ is normal we may 
consider the natural action $P_*$ of $P$ on the predual 
of $M$: $P_*(\rho)=\rho\circ P$, $\rho\in M_*$.  $P_*$ is 
a completely positive contraction on $M_*$, and its range is 
norm-closed because its adjoint $P$ has norm-closed range.   
At this point we can appeal to the elementary 
result asserting that the 
adjoint of an operator with 
norm-closed range must have weak$^*$-closed 
range.

Let $\alpha: P(M)\to P(M)$ be the restriction of $P$ to 
its range.  $\alpha$ is a unital 
surjective normal isometry that acts on a dual operator 
system.  Its 
inverse $\alpha^{-1}$ is therefore normal as well, and 
$E=\alpha^{-1}\circ P\restriction_{P(M)}$ defines a normal 
completely positive idempotent with range $P(M)$ that fixes 
the unit of $M$ and satisfies $EP=PE=P$.  Thus we have 
exhibited a unique ``polar decomposition" $P=\alpha\circ E$.  

As in the more general 
case of \cstar s discussed in the proof of 
Proposition \ref{prop6}, $P(M)$ can be made 
into a \cstar\ by introducing the multiplication 
$x\bullet y=E(xy)$, and with respect to this structure $\alpha$ 
becomes an automorphism of \cstar s.  
Moreover, since in this case $P(M)$ is weak$^*$-closed, it 
is naturally identified with the dual of the Banach 
space $P(M)_*=M_*/P(M)_\perp$, where $P(M)_\perp$ is 
the pre-annihilator of $P(M)$.  A familiar theorem of Sakai 
(\cite{sakaiC*W*}, Theorem 1.16.7) implies that the \cstar\ $P(M)$ is 
a von Neumann algebra with respect to this multiplication.  

We conclude: {\em Every quasiautomorphism $P$ of a von 
Neumann algebra $M$ has a unique representation 
$P=\alpha\circ E$ where $E:M\to M$ is a normal completely positive 
idempotent with range $P(M)$ and $\alpha$ is a $*$-automorphism of the 
natural von Neumann algebra structure of $E(M)$ associated with 
the multiplication $x\bullet y=E(xy)$, $x,y\in M$}.  
\end{rem}

The appropriate notion of local finiteness for the category of 
von Neumann algebras involves the action of normal maps on 
the predual as follows:  

\begin{defn}\label{def1}
Let $P$ be a normal completely positive map on a von Neumann 
algebra $M$ satisfying $P(\mathbf 1)=\mathbf 1$.  $P$ is said to be {\em locally finite} if for 
every normal state $\omega$ of $M$, the set of normal states 
$\{\omega\circ P^n: n=0,1,2,\dots\}$ is relatively compact 
in the norm topology of $M_*$.  
\end{defn}

Since every element of $M_*$ is a linear combination of 
normal states, we see that a map $P: M\to M$ satisfying 
the conditions of Definition \ref{def1} has the property 
the norm-closure 
of $\{\rho\circ P^n:n\geq 0\}$ is compact for 
every $\rho\in M_*$; and therefore 
$P_*(\rho)=\rho\circ P$ defines a locally finite contraction 
in $\mathcal B(M_*)$.

\begin{thm}\label{thm4}
Let $M$ be a von Neumann algebra with separable predual and let 
$P: M\to M$ be a locally finite normal completely positive 
map satisfying 
$P(\bf 1)=\bf 1$.  There is a unique 
quasiautomorphism $Q=\alpha\circ E$ of $M$ such that 
for every normal state $\rho$ of $M$, one has 
$$
\lim_{n\to\infty}\|\rho\circ P^n-\rho\circ Q^n\|=0.  
$$

The completely positive map 
$E$ is characterized as the unique 
idempotent for which there is a sequence 
$n_1<n_2<\cdots$ of positive integers such that 
\begin{equation}\label{Eq12}
\lim_{k\to\infty}\|\rho\circ P^{n_k}-\rho\circ E\|=0,\qquad 
\rho\in M_*, 
\end{equation}
$\alpha$ is 
the restriction of $P$ to the dual operator system $E(M)$, 
and $E_*(M_*)$ is the norm-closed linear span of the set 
of all maximal eigenvectors of $P_*$.  
\end{thm}

\begin{proof}
Applying Theorem \ref{thm1} to the locally finite 
contraction $P_*\in\mathcal B(M_*)$ defined by $P_*(\rho)=\rho\circ P$,  
one obtains a unique quasiautomorphism $Q_*: M_*\to M_*$ with the 
property 
$$
\lim_{n\to \infty}\|\rho\circ P^n-Q_*^n(\rho)\|=0,\qquad \rho\in M_*.  
$$
Letting $Q:M\to M$ be the adjoint of $Q_*$, one finds 
that $Q$ is a quasiautomorphism of $M$, and the 
rest follows from Theorem \ref{thm1}.  
\end{proof}

Of course, one can drop the separability hypothesis on 
the predual of $M$ at the cost of replacing the sequential 
limit (\ref{Eq12}) with an 
appropriate more general assertion.  

\section{Semigroups}\label{S:semig}

Let $X$ be a Banach space.  By a {\em contraction semigroup} we mean 
a semigroup $T=\{T_t: t\geq 0\}$ of operators on $X$ satisfying 
$\|T_t\|\leq 1$ that is strongly continuous 
in the sense that for each $x\in X$, 
the function $t\in[0,\infty) \mapsto T_tx$ 
moves continuously in the norm of $X$.  Notice that 
we have not specified that $T_0=\bf 1$, so that in general 
$T_0$ is simply an idempotent contraction.  

\begin{prop}\label{prop4}
For every contraction semigroup $T=\{T_t: t\geq0\}$ acting 
on a Banach space $X$ and every vector $x\in X$, the following are equivalent.
\begin{enumerate}
\item[(i)]The norm-closure of the orbit $\{T_tx:t\geq 0\}$ 
is compact.  
\item[(ii)]For some $s>0$, the norm closure of 
$\{x, T_sx, T_s^2x,\dots\}$ is compact.  
\end{enumerate}
\end{prop}

\begin{proof}[Sketch of Proof]
The implication (i)$\implies$(ii) is trivial.  We sketch 
the proof of (ii)$\implies$(i).  Choose $s>0$ such that the closure 
$K_x$ of $\{x,T_sx,T_s^2x,\dots\}$ is compact.  
It suffices to show that the union 
\begin{equation}\label{Eq8}
\bigcup_{0\leq r\leq s}T_rK_x 
\end{equation}
is compact, since the set (\ref{Eq8}) 
obviously contains $\{T_tx:t\geq 0\}$.  
Consider  the map of $[0,\infty)\times X$ to $X$ 
defined by 
$(t,x)\mapsto T_tx$.  This map is continuous (with 
respect to the product 
topology of $[0,\infty)\times X$ 
and the norm topology of $X$) 
because $T$ is strongly continuous.  
Since the union (\ref{Eq8}) is the range of the restriction 
of this map to the compact subspace 
$[0,s]\times K_x\subseteq [0,\infty)\times X$, 
it follows that the set (\ref{Eq8}) is compact.  
\end{proof}

\begin{defn}
A contraction semigroup $T=\{T_t:t\geq 0\}$ is said to be 
{\em locally finite} the conditions 
of Proposition \ref{prop4} are satisfied for every $x\in X$.  
\end{defn}

Notice that Proposition \ref{prop4} implies that the semigroup 
$T$ will be locally finite whenever 
$T_1$ is a locally finite contraction.

Let $T=\{T_t:t\geq0\}$ be a contraction semigroup with the 
property that $T_t$ is a quasiautomorphism of $X$ for every $t\geq 0$. 
Then $E=T_0$ is an idempotent contraction that commutes 
with $\{T_t: t\geq0\}$, and a  
straightforward argument (that we omit) shows that 
$T_tX=EX$ for every $t\geq 0$, 
that $U_t=T_t\restriction_{EX}$ is 
an automorphism of $EX$, and that we have 
the ``polar decomposition" 
\begin{equation}\label{Eq101}
T_t=U_tE,\qquad t\geq 0.  
\end{equation}
In particular, 
the most general semigroup of quasiautomorphisms $T$ is 
obtained from a semigroup $U$ of automorphisms by a direct 
sum procedure $T_t=U_t\oplus 0$, $t>0$ analogous to the one spelled 
out in Remark \ref{rem1}.  

Let $x\in X$ be a nonzero eigenvector for $T=\{T_t:t\geq0\}$.  Then 
there is a complex number $\lambda$ in the upper half-plane 
$\{z=x+iy: y\geq 0\}$ such that 
\begin{equation}\label{Eq7}
T_tx=e^{it\lambda}x,\qquad t\geq 0.  
\end{equation}
Such a $\lambda$ belongs to the point spectrum of the 
generator of $T$; we abuse notation slightly 
by writing $\sigma_p(T)$ for the set of all 
complex numbers $\lambda$ satisfying (\ref{Eq7}).  
The eigenvector $x$ is said to 
be {\em maximal} if $\|T_tx\|=\|x\|$ for $t\geq 0$; thus, 
$x$ is maximal iff $\lambda\in\mathbb R$.  
Corresponding to Theorem \ref{thm1} we have:

\begin{thm}\label{thm2}
Let $T=\{T_t: t\geq 0\}$ be a locally finite contraction semigroup acting on 
a Banach space $X$, satisfying $T_0=\bf 1$.  
There is a unique semigroup $S=\{S_t:t\geq 0\}$ 
of quasiautomorphisms such that 
\begin{equation*}
\lim_{t\to\infty}\|T_tx-S_tx\|=0,\qquad x\in X.  
\end{equation*}

Let $E$ be the idempotent $E=S_0$,   
so that 
$S_t=U_tE$ for $t\geq0$ as in (\ref{Eq101}).  
Then the 
generator of $U$ is diagonalized by 
the set of maximal eigenvectors, and 
its point spectrum is given by 
\begin{equation*}
\sigma_p(U)=\sigma_p(T)\cap \mathbb R.    
\end{equation*}
The projection $E=S_0$ is characterized as the unique 
idempotent in the set $\mathcal L$ of strong limit points 
of $\{T_t:t\geq0\}$
\begin{equation*}
\mathcal L=\bigcap_{\alpha>0} \left\{T_t: t\geq \alpha\right\}^{-{\rm strong}}, 
\end{equation*}
and $U_t$ is the restriction of $T_t$ to $EX$,  $S_t=T_tE$, $t\geq0$.  
\end{thm}

\begin{proof}[Sketch of Proof]  
The argument is merely a variation of the proof 
of Theorem \ref{thm1}, requiring little more than a change of 
notation.  Let $M$ be the closed linear span of the set 
of all maximal eigenvectors for $T$ and let $N$ be 
the space of all asymptotically null vectors 
$$
N=\{x\in X: \lim_{t\to\infty}\|T_tx\|=0\}.  
$$
In order to show that $X=N\dotplus M$, 
one first shows that the restriction $\tilde T$ of 
the semigroup $T$ to 
$M$ has the property that the identity operator belongs to the strong 
closure of $\{\tilde T_t: t\geq \alpha\}$ for every $\alpha>0$ 
by the same method used in the proof 
of Theorem \ref{thm1}.  From that, along with an appropriate 
variation of Lemma \ref{lem2} for semigroups, it follows 
that $\tilde T$ is a 
semigroup of invertible isometries with the property that 
for every $s\geq 0$, 
the inverse of $\tilde T_s$ belongs to the strong closure 
of $\{\tilde T_t: t\geq 0\}$.  Thus the 
strong closure of $\{\tilde T_t: t\geq 0\}$ contains the 
one-parameter group generated by $\{\tilde T_t$: $t\geq0\}$.   
For each vector $x\in X$ one introduces 
the set of limit points 
$$
K_\infty(x)=\bigcap_{\alpha>0}\overline{\left\{T_tx: t\geq \alpha\right\}}, 
$$
and one shows that $K_\infty(x)\subseteq M$ by the method of Theorem \ref{thm1}, 
except that now one must replace references to almost periodic functions on 
the group $\mathbb Z$ with references to almost periodic functions 
on the group $\mathbb R$.  Once one knows that 
$K_\infty(x)\subseteq M$ for every $x\in X$, the decomposition 
$X=N\dotplus M$ follows readily as in the case of single contractions.  

The remaining assertions of Theorem \ref{thm2} are straightforward.  
\end{proof}

With these general results in 
hand, one can establish a 
natural counterpart of Theorem \ref{thm4} for semigroups of 
normal completely positive maps acting on von Neumann algebras.  
We leave the explicit formulation of that result for the 
reader.  

\section{An Unstable Example}\label{S:unstable}

It is natural to ask whether {\em all} completely 
positive contractions are stable.  More precisely, can 
the key hypothesis of 
local finiteness can 
be dropped entirely 
from Theorems \ref{thm3} and \ref{thm4}, provided that 
one is willing to give up the secondary conclusions?    
For example, in the 
second group of examples described in Section \ref{S:appl1}, there 
are many homeomorphisms $h:\mathbb T\to \mathbb T$ for which  
the completely positive 
map $P$ of (\ref{Eq14}) is not locally finite and has 
no maximal eigenvectors, even though in all such cases 
$P$ is asymptotically related to a quasiautomorphism 
$Q$ as in (\ref{Eq13}).  Such examples show that 
local finiteness is not {\em necessary} for the 
stability assertion of Theorem \ref{thm3}.

We conclude by briefly describing an example of a normal 
unital completely positive 
map $P:\mathcal B(H)\to\mathcal B(H)$ 
for which the principal conclusion 
of Theorem \ref{thm4} fails in the sense that there does not 
exist a quasiautomorphism $Q:\mathcal B(H)\to \mathcal B(H)$ 
with the property 
\begin{equation}\label{Eq15}
\lim_{n\to\infty}\|\rho\circ P^n-\rho\circ Q^n\|=0, 
\end{equation}
for every normal state $\rho$ of $\mathcal B(H)$.  A significant 
feature of this example is that it very nearly satisfies 
the hypothesis of Theorem \ref{thm4} in the 
following sense: There is 
a subspace $\mathcal S\subseteq\mathcal B(H)_*$ of 
codimension one such that 
$\{\rho\circ P^n:n\geq 0\}$ is relatively norm-compact 
for every $\rho\in\mathcal S$.  

The example is based on the heat flow of 
the canonical commutation relations $\{P_t:t\geq 0\}$ 
of \cite{ccrHeat}, a semigroup of normal 
completely positive 
maps on $\mathcal B(H)$ that is {\em pure} in the sense 
that for any pair $\rho_1,\rho_2$ of normal states of 
$\mathcal B(H)$ one has 
$$
\lim_{t\to\infty}\|\rho_1\circ P_t-\rho_2\circ P_t\|=0, 
$$
while on the other hand, there is no normal 
state $\omega$ of $\mathcal B(H)$ satisfying 
$\omega\circ P_t=\omega$, $t\geq0$.  If we fix 
$t_0>0$ and set $P=P_{t_0}$, then 
\begin{equation}\label{Eq16}
\lim_{n\to\infty}\|\rho_1\circ P^n-\rho_2\circ P^n\|=0 
\end{equation}
for normal states $\rho_1, \rho_2$, and $P$ cannot leave 
any normal state $\omega$ invariant.  

The relation 
(\ref{Eq16}), together with a simple compactness 
argument, implies that for fixed $A\in\mathcal B(H)$ 
the sequence $P(A), P^2(A),P^3(A),\dots$ is asymptotically 
a scalar sequence in the sense that there is 
a sequence of complex numbers $\lambda_1,\lambda_2,\dots$ 
(which depends on $A$) such that 
$$
\lim_{n\to\infty}(P^n(A)-\lambda_n\mathbf 1)= 0
$$
in the weak operator topology.  So if there 
were a quasiautomorphism $Q$ satisfying (\ref{Eq15}), then 
the range of $Q$ would be $\mathbb C\cdot\mathbf 1$.  
Since the von Neumann algebra associated with the range 
of $Q$ in Remark \ref{rem3} is in this case $\mathbb C$ and since 
a $*$-automorphism of $\mathbb C$ 
is the identity map, 
such a quasiautomorphism $Q$ would 
simply be a normal idempotent with range 
$\mathbb C\cdot\mathbf 1$.  Thus  $Q$  
would have the form 
$Q(A)=\omega(A)\mathbf 1$, where $\omega$ 
is a normal state of $\mathcal B(H)$.  
Since $\omega\circ Q=\omega$, 
we have 
$\|\omega\circ P^n-\omega\|=\|\omega\circ P^n-\omega\circ Q^n\|$ 
for every $n\geq 1$; and 
(\ref{Eq15}) implies  that 
$\|\omega\circ P^n-\omega\circ Q^n\|\to0$ 
as $n\to\infty$.  We conclude that $\|\omega\circ P^n-\omega\|\to0$ 
as $n\to\infty$.  Therefore $\omega\circ P=\omega$, 
contradicting the second property 
of $P$ cited above.  

Finally, let $\mathcal S$ be the 
codimension one subspace of $\mathcal B(H)_*$ consisting 
of all normal linear functionals $\rho$ on $\mathcal B(H)$ 
satisfying $\rho({\bf 1})=0$.  Every $\rho\in\mathcal S$ 
can be decomposed into a sum 
$$
\rho=\lambda(\rho_1-\rho_2)+i\mu(\rho_3-\rho_4)
$$
where $\lambda$ and $\mu$ are real scalars and each $\rho_k$ is 
a normal state.  By (\ref{Eq16}), we have 
$\|\rho\circ P^n\|\to0$ as $n\to\infty$, and therefore 
$\{\rho\circ P^n:n\geq 0\}$ is relatively norm-compact 
for every $\rho$ in the subspace $\mathcal S$.  

\vfill

\bibliographystyle{alpha}
\newcommand{\noopsort}[1]{} \newcommand{\printfirst}[2]{#1}
  \newcommand{\singleletter}[1]{#1} \newcommand{\switchargs}[2]{#2#1}

\end{document}